\documentclass[12pt]{article}
 \usepackage{latexsym}
 \usepackage{amssymb}
 \usepackage{graphicx}

 \newtheorem{theorem}{Theorem}

\newtheorem{Theorem}[theorem]{Theorem}

\newtheorem{Proposition}[theorem]{Proposition}
\newtheorem{Assumption[theorem]}{Assumption}
\newtheorem{Lemma}[theorem]{Lemma}

\newcommand{\J}{{\cal J}}
\newcommand{\A}{{\cal A}}
\newcommand{\B}{{\cal B}}

\newcommand{\I}{{\cal I}}

\newcommand{\qed}{\nobreak \ifvmode \relax \else
      \ifdim\lastskip<1.5em \hskip-\lastskip
      \hskip1.5em plus0em minus0.5em \fi \nobreak
      \vrule height0.75em width0.5em depth0.25em\fi}

\def \ep{\hbox{ }\hfill$\Box$}

\begin{document}
\title{Symmetric Nonnegative Tensors and Copositive Tensors}

\author{
Liqun Qi \thanks{Email: maqilq@polyu.edu.hk. Department of Applied
Mathematics, The Hong Kong Polytechnic University, Hung Hom,
Kowloon, Hong Kong. This author's work was supported by the Hong
Kong Research Grant Council (Grant No. PolyU 501909, 502510, 502111
and 501212).}}

\date{\today} \maketitle

\begin{abstract}
\noindent We first prove two new spectral properties for symmetric
nonnegative tensors.   We prove a maximum property for the largest
H-eigenvalue of a symmetric nonnegative tensor, and establish some
bounds for this eigenvalue via row sums of that tensor. We show that
if an eigenvalue of a symmetric nonnegative tensor has a positive
H-eigenvector, then this eigenvalue is the largest H-eigenvalue of
that tensor. We also give a necessary and sufficient condition for
this. We then introduce copositive tensors. This concept extends the
concept of copositive matrices. Symmetric nonnegative tensors and
positive semi-definite tensors are examples of copositive tensors.
The diagonal elements of a copositive tensor must be nonnegative. We
show that if each sum of a diagonal element and all the negative
off-diagonal elements in the same row of a real symmetric tensor is
nonnegative, then that tensor is a copositive tensor.  Some further
properties of copositive tensors are discussed. \vspace{3mm}

\noindent {\bf Key words:}\hspace{2mm} nonnegative tensor,
copositive tensor, H-eigenvalue \vspace{3mm}

\noindent {\bf AMS subject classifications (2010):}\hspace{2mm}
15A18; 15A69
  \vspace{3mm}

\end{abstract}


\section{Introduction}
\hspace{4mm} Eigenvalues of higher-order tensors were introduced in
\cite{Qi, Li} in 2005.   Since then, many research works have been
done in spectral theory of tensors.   In particular, the theory and
algorithms for eigenvalues of nonnegative tensors are well developed
\cite{CPZ, CPZ1, FGH, HHQ, LZI, NQZ, YY, YY1, ZQ, ZQX}.

In this paper, we prove two new spectral properties for symmetric
nonnegative tensors.    Then we introduce copositive tensors,
establish some necessary conditions and some sufficient conditions
for a real symmetric tensor to be a copositive tensor, and discuss
some further properties of copositive tensors.

Some preliminary concepts and results are given in the next section.

We prove a maximum property of the largest H-eigenvalue of a
symmetric nonnegative tensor in Section 3. Based upon this, we
establish some bounds for the largest eigenvalue of a symmetric
nonnegative tensor via row sums of that tensor.

In Section 4, we show that a symmetric nonnegative tensor has at
most one H$^{++}$-eigenvalue, i.e., an H-eigenvalue with a positive
H-eigenvector, and if an eigenvalue of a symmetric nonnegative
tensor has a positive H-eigenvector, then that eigenvalue must equal
to the largest eigenvalue of that tensor.   We also give a necessary
and sufficient condition for the existence of such an
H$^{++}$-eigenvalue.

In Section 5, we introduce copositive tensors and strictly
copositive tensors. These two concepts extend the concepts of
copositive matrices and strictly copositive matrices. Symmetric
nonnegative tensors and positive semi-definite tensors are examples
of copositive tensors. The diagonal elements of a copositive tensor
must be nonnegative. We show that if each sum of a diagonal element
and all the negative off-diagonal elements in the same row of a real
symmetric tensor is nonnegative, then that tensor is a copositive
tensor.

Some further properties of copositive tensors are discussed in
Section 6.   We show that if a copositive tensor has an
H$^+$-eigenvalue, i.e., an H-eigenvalue with a nonnegative
H-eigenvector, then that H$^+$-eigenvalue must be nonnegative. The
sets of copositive tensors and strictly copositive tensors form two
convex cones: the copositive tensor cone and the strictly copositive
tensor cone.   We show that the latter is exactly the interior of
the former.   We also introduced completely positive tensors.  The
copositive tensor cone is the dual cone of the completely positive
tensor cone.   If the completely positive tensor cone is closed,
then these two cones are dual to each other.

Some final remarks are made in Section 7.

Denote by $e$ the all $1$ $n$-dimensional vector, $e_j = 1$ for
$j=1, \cdots, n$.   Denote by $e^{(i)}$ the $i$th unit vector in
$\Re^n$, i.e., $e^{(i)}_j = 1$ if $i=j$ and $e^{(i)}_j = 0$ if $i
\not =j$, for $i, j = 1, \cdots, n$. Denote the set of all
nonnegative vectors in $\Re^n$ by $\Re^n_+$ and the set of all
positive vectors in $\Re^n$ by $\Re^n_{++}$.  If both $\A =
\left(a_{i_1 \cdots i_k}\right)$ and $\B = \left(b_{i_1 \cdots
i_k}\right)$ are real $k$th order $n$-dimensional tensors, and
$b_{i_1 \cdots i_k} \le a_{i_1 \cdots i_k}$ for $i_1, \cdots, j_k =
1, \cdots, n$, then we denote $\B \le \A$.  We use $\J$ to denote
the $k$th order $n$-dimensional tensor with all of its elements
being $1$.  We use $\I$ to denote the $k$th order $n$-dimensional
diagonal tensor with all of its diagonal elements being $1$.

\section{Preliminaries}

Let $\A = \left(a_{i_1 \cdots i_k}\right)$ be a real $k$th order
$n$-dimensional tensor, and $x \in C^n$. Then
$$\A x^k = \sum_{i_1, \cdots, i_k = 1}^n a_{i_1 \cdots i_k}x_{i_1}
\cdots x_{i_k},$$ and $\A x^{k-1}$ is a vector in $C^n$, with its
$i$th component defined by
$$\left(\A x^{k-1}\right)_i = \sum_{i_2, \cdots, i_k = 1}^n a_{i i_1 \cdots i_k}x_{i_2}
\cdots x_{i_k}.$$ Let $r$ be a positive integer.  Then $x^{[r]}$ is
a vector in $C^n$, with its $i$th component defined by $x_i^r$. We
say that $\A$ is symmetric if its entries $a_{i_1, \cdots, i_k}$ are
invariant for any permutation of the indices.

If $x \in C^n$, $x \not = 0$, $\lambda \in C$,  $x$ and $\lambda$
satisfy
\begin{equation} \label{eig}
\A x^{k-1} = \lambda x^{[k-1]},
\end{equation}
then we call $\lambda$ an {\bf eigenvalue} of $\A$, and $x$ its
corresponding {\bf eigenvector}.   By (\ref{eig}), if $\lambda$ is
an eigenvalue of $\A$ and $x$ is its corresponding eigenvector, then
$$\lambda = {(\A x^{k-1})_j \over x_j^{k-1}},$$
for some $j$ with $x_j \not = 0$.   In particular, if $x$ is real,
then $\lambda$ is also real. In this case, we say that $\lambda$ is
an {\bf H-eigenvalue} of $\A$ and $x$ is its corresponding {\bf
H-eigenvector}.    If $x \in \Re^n_+$, then we say that $\lambda$ is
an {\bf H$^+$-eigenvalue} of $\A$.  If $x \in \Re^n_{++}$, then we
say that $\lambda$ is an {\bf H$^{++}$-eigenvalue} of $\A$.  The
largest modulus of the eigenvalues of $\A$ is called the {\bf
spectral radius} of $\A$, denoted by $\rho(\A)$.

By \cite{CPZ}, $\A$ is called {\bf reducible} if there exists a
proper nonempty subset $I$ of $\{ 1, \cdots, n \}$ such that
$$a_{i_1\cdots i_k} = 0, \ \ \forall i_1 \in I, \ \ \forall i_2,
\cdots , i_k \not \in I.$$ If $\A$ is not reducible, then we say
that $\A$ is {\bf irreducible}.

Let $\A = (a_{i_1, \cdots, i_k})$ be a $k$th order $n$-dimensional
symmetric nonnegative tensor.   Construct a graph $\hat G(\A) =
(\hat V, \hat E)$, where $\hat V = \cup_{j=1}^n V_j, V_j$ is a copy
of $\{ 1,\cdots, n \}$, for $j = 1, \cdots, n$. Assume that $i_j \in
V_j, i_l \in V_l, j \not = l$.   The edge $(i_j, i_l) \in \hat E$ if
and only if $a_{i_1, \cdots, i_k} \not = 0$ for some $k-2$ indices
$\{ i_1, \cdots, i_k \} \setminus \{ i_j, i_l \}$.   The tensor $\A$
is called {\bf weakly irreducible} if $\hat G(\A)$ is connected.  As
observed in \cite{FGH}, an irreducible symmetric nonnegative tensor
is always weakly irreducible.

Suppose that $\A = (a_{i_1, \cdots, i_k})$ is a $k$th order
$n$-dimensional real tensor.   If all the off-diagonal entries
$a_{i_1\cdots i_k}, (i_1, \cdots, i_k) \not = (i_1, \cdots, i_1)$
are nonnegative, then $\A$ is called an {\bf essentially nonnegative
tensor} \cite{ZQL}. If all the entries $a_{i_1\cdots i_k}$ are
nonnegative, $\A$ is called a {\bf nonnegative tensor}.   We now
summarize the Perron-Frobenius theorem for nonnegative tensors,
established in \cite{CPZ, FGH, YY}. With the new definitions of
H$^+$-eigenvalues and H$^{++}$-eigenvalues, this theorem can be
stated concisely.

\begin{theorem} {\bf (The Perron-Frobenius Theorem for Nonnegative
Tensors)} \label{t1}

(a). {\bf (Yang and Yang 2010)} If $\A$ is a nonnegative tensor of
order $k$ and dimension $n$, then $\rho(\A)$ is an H$^+$-eigenvalue
of $\A$.

(b). {\bf (Friedland, Gaubert and Han 2011)}  If furthermore $\A$ is
symmetric and weakly irreducible, then $\rho(\A)$ is the unique
H$^{++}$-eigenvalue of $\A$, with the unique eigenvector $x \in
\Re^n_{++}$, up to a positive scaling constant.

(c). {\bf (Chang, Pearson and Zhang 2008)} If moreover $\A$ is
irreducible, then $\rho(\A)$ is the unique H$^+$-eigenvalue of $\A$.

(d). {\bf (Yang and Yang 2010)} If $\A$ and $\B$ are two nonnegative
tensor of order $k$ and dimension $n$, and $\B \le \A$, then
$\rho(B) \le \rho(\A)$.

\end{theorem}

Thus, for a nonnegative tensor $\A$, its spectral radius is its
largest H-eigenvalue.

\smallskip

As observed in \cite{ZQL}, a tensor $\A$ is an essentially
nonnegative tensor, if and only if there are a nonnegative tensor
$\B$ and a real number $c$, such that $\A = \B + c\I$.

\section{A Maximum Property of the Largest H-Eigenvalue of a Symmetric Nonnegative Tensor}

Suppose that $\A = (a_{i_1\cdots i_k})$ is a $k$th order
$n$-dimensional real symmetric tensor, with $k \ge 2$. Denote its
largest H-eigenvalue by $\lambda_{\max}(\A)$. When $k$ is even, by
\cite{Qi}, we know that
$$\lambda_{\max}(\A) = \max \{ \A x^k : x \in \Re^n, \sum_{i=1}^n x_i^k = 1
\}.$$

In this section, we first prove the following theorem, which holds
whenever $k$ is even or odd.

\begin{theorem} {\bf (A Maximum Property of The Largest H-Eigenvalue of a Symmetric Nonnegative Tensor)}
\label{t2} Suppose that $\A = (a_{i_1\cdots i_k})$ is a $k$th order
$n$-dimensional symmetric nonnegative tensor, with $k \ge 2$.
 Then we have
 \begin{equation} \label{e2.1}
\lambda_{\max}(\A) = \max \{ \A x^k : x \in \Re^n_+, \sum_{i=1}^n
x_i^k = 1 \}.
\end{equation}
\end{theorem}

\noindent {\bf Proof.} We now prove (\ref{e2.1}).  Assume that $\A$
is a symmetric nonnegative tensor.  By Theorem \ref{t1}, we have
$$\lambda_{\max}(\A) = \max \{ \lambda : \A x^{k-1} = \lambda x^{[k-1]}, x \in
\Re^n_+ \} = \max \{ \lambda : \A x^{k-1} = \lambda x^{[k-1]}, x \in
\Re^n_+, \sum_{i=1}^n x_i^k = 1 \}$$
$$ = \max \{ \A x^k : \A x^{k-1} = \lambda x^{[k-1]}, x \in
\Re^n_+, \sum_{i=1}^n x_i^k = 1 \} \le \max \{ \A x^k : x \in
\Re^n_+, \sum_{i=1}^n x_i^k = 1 \}.$$

On the other hand, assume that $x^*$ is an optimal solution of $\max
\{ \A x^k : x \in \Re^n_+, \sum_{i=1}^n x_i^k = 1 \}$.    By
optimization theory, there is a Lagrangian multiplier $\lambda$ and
a nonempty subset $I$ of $\{ 1, \cdots, n \}$ such that for $i \in
\{ 1, \cdots, n \} \setminus I$, $x_i^* = 0$, and for $i \in I$,
$$\left(\A (x_i^*)^{k-1}\right)_i = \lambda (x_i^*)^{k-1}.$$
Multiplying the above equalities by $x_i^*$ and summing up them, we
have
$$\lambda = \A (x^*)^k = \max \{ \A x^k : x \in \Re^n_+, \sum_{i=1}^n x_i^k = 1
\}.$$ Construct a $k$th order $n$-dimensional symmetric nonnegative
tensor $\B = (b_{i_1\cdots i_k})$ such that $b_{i_1\cdots i_k} =
a_{i_1\cdots i_k}$ if $i_1, \cdots, i_k \in I$, and $b_{i_1\cdots
i_k} = 0$ otherwise.  Then we see that $\lambda$ is an H-eigenvalue
of $\B$ with an H-eigenvector $x^*$.   Then we see that $ \B  \le
\A$.   By Theorem \ref{t1} (d), we have
$$\lambda \le \rho(\B) \le \rho(\A) = \lambda_{\max}(\A).$$
Combining these together, we have (\ref{e2.1}). \ep

The adjacency tensor of a uniform hypergraph is a nonnegative tensor
\cite{CD}.   The signless Laplacian tensor of a uniform hypergraph,
introduced in \cite{Qi1}, is also a nonnegative tensor.  Cooper and
Dutle \cite{CD} established (\ref{e2.1}) for the adjacency tensor of
a connected uniform hypergraph.   Qi \cite{Qi1} established
(\ref{e2.1}) for the adjacency tensor and the signless Laplacian
tensor of a general uniform hypergraph.   Zhang \cite{Za} pointed
out that (\ref{e2.1}) holds for a weakly irreducible symmetric
nonnegative tensor.   Here, we established (\ref{e2.1}) for a
general symmetric nonnegative tensor.

With Theorem \ref{t2}, we may establish some lower bounds for
$\rho(\A)$.    We define the $i$th {\bf row sum} of a $k$th order
$n$-dimensional tensor $\A = (a_{i_1, \cdots, i_k})$ as
$$R_i(\A) = \sum_{i_2,\cdots, i_k = 1}^n a_{ii_2\cdots i_k},$$
and denote the largest, the smallest and the average row sums of
$\A$ by
$$R_{\max}(\A) = \max_{i = 1, \cdots, n} R_i(\A), \ R_{\min}(\A) = \min_{i = 1, \cdots, n} R_i(\A), \ {\rm and}\  \bar R(\A) = {1 \over n}\sum_{i = 1}^n R_i(\A),$$
respectively.   We also denote the largest, the smallest and the
average diagonal element of $\A$ by
$$d_{\max}(\A) = \max_{i = 1, \cdots, n} a_{i\cdots i}, \ d_{\min}(\A) = \min_{i = 1, \cdots, n} a_{i\cdots i},
\ \bar d(\A) = {1 \over n}\sum_{i = 1}^n a_{i\cdots i},$$
respectively.

\begin{theorem} {\bf (Bounds The Largest H-Eigenvalue of a Nonnegative Tensor)}\label{t3}
Suppose that $\A = (a_{i_1\cdots i_k})$ is a $k$th order
$n$-dimensional nonnegative tensor, with $k \ge 2$.
 Then we have
 \begin{equation} \label{e2.2}
\lambda_{\max}(\A) \le R_{\max}(\A).
\end{equation}
If furthermore $\A$ is symmetric, then we have
 \begin{equation} \label{e2.3}
\lambda_{\max}(\A) \ge \max \{ \bar R(\A), d_{\max}(\A) \}.
\end{equation}
\end{theorem}
\noindent {\bf Proof.} By Theorem \ref{t1}, $\A$ has a nonnegative
H-eigenvector $x$.  Let $x_j = \max_{i=1, \cdots, n} x_i$.   Then
$x_j > 0$.   We have
$$\sum_{i_2, \cdots, i_k=1}^n a_{ji_2\cdots i_k}x_{i_2}\cdots
x_{i_k} = \lambda_{\max}(\A)x_j^{k-1},$$ i.e.,
$$\lambda_{\max}(\A) = \sum_{i_2, \cdots, i_k=1}^n a_{ji_2\cdots i_k}{x_{i_2}\over x_j}\cdots
{x_{i_k} \over x_j} \le R_j(\A) \le R_{\max}(\A).$$ This proves
(\ref{e2.2}).

Now assume that $\A$ is symmetric.   Let $y = {e \over (n)^{1 \over
k}}$.   By Theorem \ref{t2}, we have
$$\lambda_{\max}(\A) \ge \A y^k = {1 \over n}\sum_{i_1, \cdots, i_k
= 1}^n a_{i_1\cdots i_k} = \bar R(\A).$$ Assume that $a_{j\cdots j}
= d_{\max}(\A)$. Let $y = e^{(j)}$.  By Theorem \ref{t2}, we have
$$\lambda_{\max}(\A) \ge \A y^k = a_{j\cdots j} = d_{\max}(\A).$$
Combining these two inequalities, we have (\ref{e2.3}). \ep

If we apply Theorem \ref{t3} to the adjacency tensor of a uniform
hypergraph, we may get the bounds for the largest H-eigenvalue of
that tensor, obtained by Cooper and Dutle in \cite{CD}. If we apply
Theorem \ref{t3} to the signless Laplacian tensor of a uniform
hypergraph, we may get the bounds for the largest H-eigenvalue of
that tensor, obtained by Qi in \cite{Qi1}.

For a $k$th order $n$-dimensional symmetric nonnegative tensor $\A$,
if all of its row sums are the same, then we have
$\lambda_{\max}(\A) = \bar R(\A)$.   The adjacency tensor and the
signless Laplacian tensor of a regular $k$-graph are such examples
\cite{CD, Qi1}.

\section{The H$^{++}$-Eigenvalue of a Symmetric Nonnegative Tensor}

In this section, we show that a symmetric nonnegative tensor has at
most one H$^{++}$-eigenvalue.

Suppose that $I \subset \{ 1, \cdots, n \}$.   Let $x$ be an
$n$-dimensional vector.  Then $x(I)$ is an $|I|$-dimensional vector
with its components indexed for $i \in I$, and $x(I)_i \equiv x_i$
for $i \in I$.   For a $k$th order $n$-dimensional tensor $\A =
(a_{i_1\cdots i_k})$, $\A(I)$ is a $k$th order $|I|$-dimensional
tensor with elements $a_{i_1\cdots i_k}, i_1, \cdots, i_k \in I$.

Suppose that $\A = (a_{i_1\cdots i_k})$ is a symmetric nonnegative
tensor of order $k$ and dimension $n$.  By \cite{FGH, HHQ}, there is
a partition $(I_1, \cdots, I_s)$ of $\{ 1, \cdots, n \}$, such that
$\A(I_r)$ is weakly irreducible for $r = 1, \cdots, s$, and
$a_{i_1\cdots i_k} = 0$ for all $i_1 \in I_r, i_2, \cdots, i_k \not
\in I_r, r = 1, \cdots, s$.    Furthermore, we have
$$\lambda_{\max}(\A) = \max \{ \lambda_{\max}(\A(I_r)) : r = 1, \cdots, s \}.$$

\begin{Theorem} \label{t4} {\bf (The H$^{++}$-Eigenvalue of a Symmetric Nonnegative Tensor)}
Let $\A = (a_{i_1\cdots i_k})$ be a symmetric nonnegative tensor of
order $k$ and dimension $n$.  Then $\A$ has at most one
H$^{++}$-eigenvalue. A real number $\lambda$ is an
H$^{++}$-eigenvalue of $\A$ if and only if for the above partition
$(I_1, \cdots, I_s)$, we have
\begin{equation} \label{e4.1}
\lambda = \lambda_{\max}(\A) = \lambda_{\max}(\A(I_r)), \ {\rm for}\
r = 1, \cdots, s.
\end{equation}
\end{Theorem}

\noindent {\bf Proof.}   Suppose that (\ref{e4.1}) holds.  Then by
Theorem \ref{t1} (a),  we have $x \in \Re^n_+$ such that
\begin{equation} \label{e4.2}
\sum_{i_2, \cdots, i_k \in I_r} a_{ii_2\cdots i_k}x_{i_2}\cdots
x_{i_k} = \lambda x_i^{k-1},
\end{equation}
for $i \in I_r, r = 1, \cdots, s$. By Theorem \ref{t1} (b), $x(I_r)
> 0$ for $r = 1, \cdots, s$.  Thus, $x \in \Re^n_+$.   (\ref{e4.2})
further implies that
\begin{equation} \label{e4.3}
\sum_{i_2, \cdots, i_k = 1}^n a_{ii_2\cdots i_k}x_{i_2}\cdots
x_{i_k} = \lambda x_i^{k-1},
\end{equation}
for $i = 1, \cdots, n$, i.e., $\lambda$ is an H$^{++}$-eigenvalue of
$\A$.

On the other hand, assume that $\lambda$ is an H$^{++}$-eigenvalue
of $\A$, with an H-eigenvector $x \in \Re^n_{++}$.  Then we have
(\ref{e4.3}), which implies (\ref{e4.2}).  By Theorem \ref{t1} (b),
we have $\lambda = \lambda_{\max}(\A(I_r))$ for $r = 1, \cdots, s$.
Since $\lambda_{\max}(\A) = \max \{ \lambda_{\max}(\A(I_r)) : r = 1,
\cdots, s \}$, we have (\ref{e4.1}).  \ep

\section{Copositive Tensors}

The concept of copositive matrices was introduced by Motzkin
\cite{Mo} in 1952.   It is an important concept in applied
mathematics, with applications in control theory, optimization
modeling, linear complementarity problems, graph theory and linear
evolution variational inequalities \cite{HS}. We now extend this
concept to tensors.

Suppose that  $\A = (a_{i_1\cdots i_k})$ is a real symmetric tensor
of order $k$ and dimension $n$.  We say that $\A$ is a {\bf
copositive tensor} if for any $x \in \Re^n_+$, we have $\A x^k \ge
0$.  We say that $\A$ is a {\bf strictly copositive tensor} if for
any $x \in \Re^n_+, x \not = 0$, we have $\A x^k > 0$.   Clearly, a
symmetric nonnegative tensor is a copositive tensor.  Recall
\cite{Qi} that a real symmetric tensor $\A$ of order $k$ and
dimension $n$, is called a {\rm positive semi-definite tensor}, if
for any $x \in \Re^n$, $\A x^k \ge 0$, $\A$ is called a {\rm
positive definite tensor}, if for any $x \in \Re^n, x \not = 0$, $\A
x^k > 0$.  Except the zero tensor, positive semi-definite tensors
are of even order.   Clearly, a positive semi-definite tensor is a
copositive tensor, a positive definite tensor is a strictly
copositive tensor.

\begin{Theorem} \label{t5} {\bf (Copositive Tensors)} Suppose that  $\A = (a_{i_1\cdots i_k})$ and $\B = (b_{i_1\cdots i_k})$ are two real symmetric
tensors of order $k$ and dimension $n$.   Then we have the following
conclusions.

(a). $\A$ is copositive if and only if
\begin{equation} \label{e5.1}
N_{\min}(\A) \equiv \min \{ \A x^k : x \in \Re^n_+, \sum_{i=1}^n
x_i^k = 1 \} \ge 0.
\end{equation}
$\A$ is strictly copositive if and only if
\begin{equation} \label{e5.2}
N_{\min}(\A) \equiv \min \{ \A x^k : x \in \Re^n_+, \sum_{i=1}^n
x_i^k = 1 \} > 0.
\end{equation}

(b). If $\A$ is copositive, then $d_{min}(\A) \ge 0$.  If $\A$ is
strictly copositive, then $d_{min}(\A) > 0$.

(c). Suppose that $\A \le \B$.   If $\A$ is copositive, then $\B$ is
copositive.  If $\A$ is strictly copositive, then $\B$ is strictly
copositive.
\end{Theorem}

\noindent {\bf Proof.} (a). If $\A$ is copositive, then clearly
(\ref{e5.1}) holds.   Suppose (\ref{e5.1}) holds.   For any $y \in
\Re^n_+, y \not = 0$, let
$$x = {y \over \left(\sum_{i=1}^n y_i^k\right)^{1
\over k}}.$$ Then $x \in \Re^+$, $\sum_{i=1}^n x_i^k = 1$, and
$$ \A x^k = {\A y^k \over \sum_{i=1}^n y_i^k} \ge {N_{\min}(\A) \over \sum_{i=1}^n y_i^k} \ge 0.$$
Thus, $\A$ is copositive.

Similarly, if (\ref{e5.2}) holds, we may show that $\A$ is strictly
copositive.  Suppose that (\ref{e5.2}) does not hold.   As the
feasible set of the minimization problem in (\ref{e5.2}) is compact,
the minimization problem has an optimizer $x^*$.  Then $x^* \in
\Re^n_+, x^* \not = 0$ and $\A (x^*)^k = N_{\min}(\A) \le 0$.   Thus
$\A$ cannot be strictly copositive.  This completes the proof of
(a).

(b). Assume that $d_j(\A) = d_{\min}(\A)$.  Let $y = e^{(j)}$. Then
$y \in \Re^n_+, \sum_{i=1}^n y_i^k = 1$, and $d_{\min}(\A) = \A
y^k$.  If $\A$ is copositive, then by (a),
$$d_{\min}(\A) = \A y^k \ge N_{\min}(\A) \ge 0.$$
If $\A$ is strictly copositive, then by (a),
$$d_{\min}(\A) = \A y^k \ge N_{\min}(\A) > 0.$$  These prove (b).

(c).  Suppose that $\A \le \B$.  If $\A$ is copositive, then for any
$x \in \Re^n_+$, $\B x^k \ge \A x^k \ge 0$.  This implies that $\B$
is copositive.  If $\A$ is strictly copositive, then for any $x \in
\Re^n_+, x \not = 0$, $\B x^k \ge \A x^k > 0$.  This implies that
$\B$ is strictly copositive. \ep

We now prove further a nontrivial sufficient condition for a real
symmetric tensor to be copositive.   We need to prove some lemmas
first.

\begin{Lemma} \label{l1} Suppose that $\A = (a_{i_1\cdots i_k})$ is a $k$th order
$n$-dimensional symmetric, essentially nonnegative tensor, with $k
\ge 2$.  Then we still have (\ref{e2.1}).
\end{Lemma}

\noindent {\bf Proof.}  Assume that $\A$ is a symmetric, essentially
nonnegative tensor. Then there are a symmetric nonnegative tensor
$\B$ and a real number $c$, such that $\A = \B + c\I$.  Then,
$$\lambda_{\max}(\A) = \lambda_{\max}(\B) + c = \max \{ \B x^k : x \in \Re^n_+, \sum_{i=1}^n x_i^k = 1
\}+ c$$ $$= \max \{ \B x^k + c : x \in \Re^n_+, \sum_{i=1}^n x_i^k =
1 \} = \max \{ \A x^k : x \in \Re^n_+, \sum_{i=1}^n x_i^k = 1 \}.$$
This completes the proof. \ep

We may also show that Theorem \ref{t3} also holds for symmetric,
essentially nonnegative tensors, and Theorem \ref{t4} also holds for
essentially nonnegative tensors.   We do not go to the details.

Suppose that $\A = (a_{i_1, \cdots, i_k})$ is a $k$th order
$n$-dimensional real tensor.   If all the off-diagonal entries
$a_{i_1\cdots i_k}, (i_1, \cdots, i_k) \not = (i_1, \cdots, i_1)$
are nonpositive, then $\A$ is called an {\bf essentially nonpositive
tensor}.

For a $k$th order $n$-dimensional real tensor $\A$, denote its
smallest H-eigenvalues by $\lambda_{\min}(\A)$.

\begin{Lemma} \label{l2} Suppose that $\A = (a_{i_1\cdots i_k})$ is a $k$th order
$n$-dimensional symmetric, essentially nonpositive tensor, with $k
\ge 2$. Then we have
 \begin{equation} \label{e5.3}
\lambda_{\min}(\A) = \min \{ \A x^k : x \in \Re^n_+, \sum_{i=1}^n
x_i^k = 1 \}.
\end{equation}
\end{Lemma}
\noindent {\bf Proof.}  Assume that $\A$ is a symmetric, essentially
nonpositive tensor.   Let $\B = - \A$.  Then $\B$ is a symmetric,
essentially nonnegative tensor. Then,
$$\lambda_{\min}(\A) = - \lambda_{\max}(\B) = - \max \{ \B x^k : x \in \Re^n_+, \sum_{i=1}^n x_i^k = 1
\} = \min \{ \A x^k : x \in \Re^n_+, \sum_{i=1}^n x_i^k = 1 \}.$$
This completes the proof. \ep

\begin{Lemma} \label{l2.1} Suppose that $\A = (a_{i_1\cdots i_k})$ is a $k$th order
$n$-dimensional essentially nonpositive tensor, with $k \ge 2$. Then
we have
$$\lambda_{\min}(\A) \ge R_{\min}(\A).$$
If furthermore $\A$ is symmetric, then we have
$$\lambda_{\min}(\A) \le \min \{ \bar R(\A), d_{\min}(\A) \}.$$
\end{Lemma}

The proof of this lemma is similar to the proof of Theorem \ref{t3}.
By Lemma \ref{l2.1} and Theorem \ref{t5} (a), we have the following
lemma.

\begin{Lemma} \label{l3} Suppose that $\A = (a_{i_1\cdots i_k})$ is a $k$th order
$n$-dimensional symmetric, essentially nonpositive tensor, with $k
\ge 2$. If $R_{\min}(\A) \ge 0$, for $i = 1, \cdots, n$, then $\A$
is copositive. If $R_{\min}(\A) > 0$, for $i = 1, \cdots, n$, then
$\A$ is strictly copositive.
\end{Lemma}

Finally, we may prove the following theorem.

\begin{Theorem} \label{t6} Suppose that $\B = (b_{i_1\cdots i_k})$ is a $k$th order
$n$-dimensional real symmetric tensor, with $k \ge 2$. If
 \begin{equation} \label{e5.6}
b_{i\cdots i} + \sum \{ b_{ii_2\cdots i_k} : b_{ii_2\cdots i_k} < 0,
(i, i_2, \cdots, i_k) \not = (i, \cdots, i) \} \ge 0,
\end{equation}
for $i = 1, \cdots, n$, then $\B$ is copositive. If
 \begin{equation} \label{e5.7}
b_{i\cdots i} + \sum \{ b_{ii_2\cdots i_k} : b_{ii_2\cdots i_k} < 0,
(i, i_2, \cdots, i_k) \not = (i, \cdots, i) \} > 0,
\end{equation}
for $i = 1, \cdots, n$, then $\B$ is strictly copositive.
\end{Theorem}
\noindent {\bf Proof.} Construct a $k$th order $n$-dimensional real
symmetric tensor $\A = (a_{i_1\cdots i_k})$ by $a_{i_1\cdots i_k} =
0$ if $b_{i_1\cdots i_k} > 0$ and $(i_1, \cdots, i_k) \not = (i_1,
\cdots, i_1)$, and $a_{i_1\cdots i_k} = b_{i_1\cdots i_k}$
otherwise. Then $\A$ is symmetric and essentially nonpositive, and
$\A \le \B$. Now the conclusions follow from Theorem \ref{t5} (c)
and Lemma \ref{l3}. \ep

We may call a real symmetric tensor satisfying (\ref{e5.6}) a {\bf
nonnegative diagonal dominated tensor}, and a real symmetric tensor
satisfying (\ref{e5.7}) a {\bf positive diagonal dominated tensor}.
We see that a symmetric nonnegative tensor is a nonnegative diagonal
dominated tensor.  The Laplacian tensor of a uniform hypergraph,
introduced in \cite{Qi1}, is a nonnegative diagonal dominated
tensor.  Thus, the adjacency tensor, the Laplacian tensor and the
signless Laplacian tensor of a uniform hypergraph are examples of
copositive tensors.

\section{Further Properties of Copositive Tensors}

The proofs of many properties of copositive matrices may not be
extended to copositive tensors directly.   This leaves some puzzles:
do such properties of copositive matrices still hold for copositive
tensors?    The situation is in particular odd when the order $k$ is
odd.

The first question is: Does a copositive tensor have an
H-eigenvalue?   When the order $k$ is even, by \cite{Qi}, a real
symmetric $k$th order $n$-dimensional tensor $\A$ always has an
H-eigenvalue.  Thus, we may ask

{\bf Question 1}. When the order $k$ is odd, does a copositive
tensor $\A$ always have an H-eigenvalue?

{\bf Question 2}. When the order $k \ge 3$, if a copositive tensor
$\A$ has at least one H-eigenvalue, does it always have a
nonnegative H-eigenvalue?

For a copositive matrix $A$, Haynsworth and Hoffman \cite{HH} showed
that its largest eigenvalue $\lambda$ satisfies that $\lambda \ge |
\mu|$, where $\mu$ is any other eigenvalue of $A$.   It is not clear
if this is still true for copositive tensors.

{\bf Question 3}. When the order $k \ge 3$, if a copositive tensor
$\A$ has a nonnegative H-eigenvalue $\lambda$, does it satisfy
$\lambda \ge | \mu|$, where $\mu$ is any other H-eigenvalue of $\A$?

If $\A$ is a symmetric nonnegative tensor, then the answers to the
above three questions are all ``yes''.  If $\A$ is a symmetric,
essentially nonpositive tensor, and $R_{\min}(\A) \ge 0$, then the
answers to the above three questions are also all ``yes''.

When $k$ is even, if all the H-eigenvalues of a real symmetric
tensor are nonnegative, then that tensor is positive semi-definite,
thus copositive.   The situation becomes again puzzled when $k$ is
even.

{\bf Question 4}.  Suppose that the order $k$ is odd and all the
H-eigenvalues of a real symmetric tensor are nonnegative.   Is that
tensor copositive?

No matter such basic questions remain open, we may derive some
further properties of a copositive tensor.

\begin{Proposition} \label{p1}
If a copositive tensor $\A$ has an H$^+$-eigenvalue $\lambda$, then
$\lambda \ge 0$.
\end{Proposition}

\noindent {\bf Proof.}  Then we have $\A x^{k-1} = \lambda
x^{[k-1]}$ with $x \in \Re^n_+, x \not = 0$.  We have $\lambda = {\A
x^k \over \sum_{i=1}^n x_i^k} \ge 0$. \ep

We now extend one theorem of V\"{a}liaho, Theorem 3.2 of \cite{Va},
to copositive tensors.

\begin{Proposition} \label{p2} Suppose that $\A$ is a $k$th order $n$-dimensional
copositive tensor.  Then $x \in \Re^n_+$ and $\A x^k = 0$ imply that
$\A x^{k-1} \ge 0$.
\end{Proposition}

\noindent {\bf Proof.}  Consider $f(x) = \A x^k$.  If $\A x^k = 0$
for some $x \in \Re^n_+$, then for $ t> 0$ and $i = 1, \cdots, n$,
$\A (x+te^{(i)})^k \ge 0$ as $\A$ is copositive and $x +te^{(i)} \in
\Re^n_+$.   This implies that $f'(x) = k\A x^{k-1} \ge 0$.   \ep

However, it is not clear if the next theorem of V\"{a}liaho, Theorem
3.3 of \cite{Va}, can be extended to copositive tensors or not. This
leaves another puzzle.

It is easy to see that if $\A$ and $\B$ are two (strictly)
copositive tensors of the same order and dimension, then $\A + \B$
is also a (strictly) copositive tensor, and if $\A$ is a (strictly)
copositive tensor and $\alpha$ is a positive number, then $\alpha
\A$ is also a (strictly) copositive tensor. Then all copositive
tensors of order $k$ and dimension $n$ form a convex cone. We denote
it by $C_{k, n}$.  Similarly,  all strictly copositive tensors of
order $k$ and dimension $n$ form a convex cone. We denote it by
$SC_{k, n}$.   Similarly, we have the positive semi-definite tensor
cone of order $k$ and dimension $n$, denoted by $PSD_{k, n}$, and
the nonnegative diagonal dominated tensor cone of order $k$ and
dimension $n$, denoted by $NDD_{k, n}$, etc.   When $k$ is odd,
$NDD_{k, n}$ is a subcone of $C_{k, n}$. When $k$ is even, $NDD_{k,
n}$ and $PSD_{k, n}$ are two subcones of $C_{k, n}$.

Even when $k$ is odd, a copositive tensor may not be a nonnegative
diagonal dominated tensor.   For example, let $k=n=3$, $a_{113}=
a_{131} = a_{311} = a_{223} = a_{232}=a_{322}=2$, $a_{123} = a_{132}
= a_{213} = a_{231}=a_{312}=a_{321}=-1$, and the other elements of
$\A$ be zero.  Then $\A x^3 = 6(x_1^2+ x_2^2 - x_1x_2)x_3 \ge 0$ for
any $x \in \Re^3_+$, i.e., $\A$ is a copositive tensor. But $\A$ is
not a nonnegative diagonal dominated tensor, as all diagonal
elements of $\A$ are zero, but there are negative off-diagonal
elements.

\begin{Proposition} \label{p3}
$SC_{k, n}$ is exactly the interior cone of $C_{k, n}$.
\end{Proposition}

\noindent {\bf Proof.}  Denote $B_{k, n}$ as the set of all $k$th
order $n$-dimensional real symmetric tensors whose Frobenius norms
are $1$.    Suppose that $\A \in SC_{k, n}$.  Let $\A(t, \B) = \A +
t\B$, where $\B \in B_{k, n}$.   Let $\delta$ be a positive number,
$0 \le t \le \delta$.    Then by (\ref{e5.2}),
$$\left| N_{\min}(\A(t, \B)) - N_{\min}(\A)\right| \le c\delta,$$
where $c$ is a certain norm ratio constant.   Thus, we have some
$\delta > 0$, such that for all $\B \in B_{k, n}$ and $0 \le t \le
\delta$, $\A(t, \B) \in SC_{k, n}$.  This shows that $SC_{k, n}$ is
in the interior of $C_{k, n}$.

On the other hand, suppose that $\A \in C_{k, n} \setminus SC_{k,
n}$.   By Theorem \ref{t5}, $N_{\min}(\A) = 0$.  Then there is a $y
\in \Re^n_+$, such that $\sum_{i=1}^n y_i^k = 1$ and $\A y^k = 0$.
Let $\A(t) = \A - ty^k$.  Then we see that $N_{\min}(\A(t)) < 0$ for
all $t > 0$.  Thus, $\A$ is not in the interior of $C_{k, n}$. This
completes our proof.  \ep

It is well-known that the copositive matrix cone and the completely
positive matrix cone are dual to each other \cite{BS, HS, Xu}. This
was established by Hall and Newman \cite{HN}.   We may consider this
issue in the tensor case. Let $y \in \Re^n_+$.  Then we may regard
$y^k$ as a rank-one $k$th order $n$-dimensional completely positive
tensor $y^k = (y_{i_1}\cdots y_{i_k})$.  We call a $k$th order
$n$-dimensional tensor $\A$ a {\bf completely positive tensor} if
there are $y^{(1)}, \cdots, y^{(r)} \in \Re^n_+$ such that
$$\A = \sum_{i=1}^r \left( y^{(i)}\right)^k.$$
The smallest value of $r$ to make the above expression hold is
called the {\bf CP-rank} of $\A$.   Clearly, a completely positive
tensor is a symmetric nonnegative tensor, and all the $k$th order
$n$-dimensional completely positive tensors form a convex cone, the
completely positive tensor cone, denoted as $CP_{k, n}$.   For two
$k$th order $n$-dimensional real symmetric tensors $\A =
(a_{i_1\cdots i_k})$ and $\B = (b_{i_1\cdots i_k})$, denote its
inner product as
$$\langle \A, \B \rangle = \sum_{i_1, \cdots, i_k=1}^n a_{i_1\cdots i_k}b_{i_1\cdots i_k}.$$
Denote the space of all $k$th order $n$-dimensional real symmetric
tensors as $S_{k, n}$.   For a convex cone $K$ in $S_{k, n}$, its
dual cone is defined as
$$K^* = \left\{ \B \in S_{k, n} : \langle \A, \B \rangle \ge 0
\right\}.$$ We have $K^{**} = {\rm cl} K$.   If $K$ is closed, then
we have $K^{**} = K$.   By the definition of copositive tensors, we
have $C_{k, n} = CP_{k, n}^*$.   Then we have $C_{k, n}^* = {\rm
cl}CP_{k, n}$.  If $CP_{k, n}$ is closed, then we have $C_{k, n}^* =
CP_{k, n}$.  We leave this as a future research topic.

\section{Final Remarks}

In Sections 3 and 4, we established two new spectral properties of
symmetric nonnegative tensors.  This shows that there are still
unexplored topics of the spectral theory of symmetric nonnegative
tensors. In Section 5, we introduced copositive tensors and strictly
copositive tensors. Symmetric nonnegative tensors and positive
semi-definite tensors are copositive tensors.   Beside some simple
properties of copositive tensors and strictly copositive tensors, we
show that nonnegative diagonal dominated tensors are copositive
tensors, and positive diagonal dominated tensors are strictly
copositive tensors. Section 6 shows that there are many puzzles
unsolved for copositive tensors and strictly copositive tensors.
Hence, this paper is only a starting point for studying copositive
tensors and strictly copositive tensors.



\end{document}